\documentclass[11pt,fullpage]{article}

\usepackage{amsmath}
\usepackage{graphicx}
\usepackage{amsfonts}
\usepackage{amssymb}
    \usepackage{amsthm}

\usepackage[normalem]{ulem}
\usepackage{scrextend}
\usepackage{endnotes}
\usepackage[usenames]{color}
\numberwithin{equation}{section}
\usepackage{natbib}
 \bibpunct[, ]{(}{)}{,}{a}{}{,}%

\newtheorem{theorem}{Theorem}

\renewcommand{\baselinestretch}{1.2}

\newcommand{\e}{\mathbb{E}}

\newcommand{\var}{\mathbb{V}\mbox{ar}}

\newcommand{\Pro}{\mathbb{P}}
\newcommand{\cov}{\mathbb{C}\mbox{ov}}
\newcommand{\cL}{\mathcal{L}}
\newcommand{\cE}{\mathcal{E}}

\title{\LARGE Gaussian Limits for Scheduled Traffic with Super-Heavy Tailed Perturbations    \vspace{0.5cm}
}

\author{{\large   Victor F. Araman \hspace{0.5cm}  Peter W. Glynn}
\thanks{The first author is with the  Olayan School of Business, American University of Beirut, Beirut, va03@aub.edu.lb. The second
author is with the Department of Management Science and Engineering at Stanford University, Stanford, CA, 74305, glynn@stanford.edu.}}

\begin{document}
    \maketitle \thispagestyle{empty}
\date{}
\renewcommand{\thefootnote}{\fnsymbol{footnote}}
\setcounter{footnote}{1}
\renewcommand{\baselinestretch}{1.2} {\small\normalsize}

\begin{abstract}    
A scheduled arrival model is one in which the $j^{th}$ customer is scheduled to arrive at time $jh$ but the customer actually arrives at time $jh + \xi_j$, where the $\xi_j$’s are independent and identically distributed. It has previously been shown that the arrival counting process for scheduled traffic obeys a functional central limit theorem (FCLT) with fractional Brownian motion (fBM) with Hurst parameter $H \in (0,1/2)$ when the $\xi_j$’s have a Pareto-like tail with tail exponent lying in $(0,1)$. Such limit processes exhibit less variability than Brownian motion, because the scheduling feature induces negative correlations in the arrival process. In this paper, we show that when the tail of the $\xi_j$’s have a super-heavy tail, the FCLT limit process is Brownian motion (i.e. $H=1/2$), so that the heaviness of the tails eliminates any remaining negative correlations and generates a limit process with independent increments. We further study the case when the $\xi_j$’s have a Cauchy-like tail, and show that the limit process in this setting is a fBM with $H=0$. So, this paper shows that the entire range of fBMs with $H \in [0,1/2]$ are possible as limits of scheduled traffic. 
\end{abstract}

\maketitle

\section{Introduction}
An important element in describing a queue is the process that characterizes the arrival of customers to the system. In this paper, we consider an arrival model known as scheduled traffic that is quite natural from a modeling viewpoint, especially within systems in which customer arrivals are governed by an appointments-based system. Specifically, a scheduled traffic model is one which the $n^{th}$ customer is scheduled to arrival at time $n h$ (with $h>0$), but her actual arrival occurs at time $nh+\xi_n.$ We call the random variable (rv) $\xi_n$ the perturbation associated with customer $n$'s arrival time. In many applications, it seems reasonable to assume the $\xi_n$'s can be modeled as a sequence of independent and identically distributed (iid) random variables (rvs). Relative to renewal traffic, which is often poorly motivated as an arrival model, scheduled traffic seems well suited to many applied domains. 

Scheduled traffic was first analyzed by \citet{Winsten59}, and is also discussed in the early queueing book by \citet{CoxSmith1961}; see their discussion of regular arrivals with unpunctuality. Most of the subsequent work has been restricted to bounded perturbations. The waiting time distribution was a focus of \citet{Winsten59}, \citet{Mercer60}, \citet{Mercer60}, \citet{Loynes62} and \citet{Mercer73}. However, their analyses did not lend themselves to direct quantitative computation. \citet{Kingman62} obtained a heavy-traffic result for single server queues with general arrival processes, a special case of which is scheduled traffic with positive finite mean perturbations (and which is noted there as being a model of special importance). 

\citet{Chen97} used scheduled traffic with deterministic service times to model aircraft landings and looked at the stability of the corresponding single server queue in the critically loaded regime when the perturbations are bounded and service times are deterministic. \citet{AramanGlynn23} have recently extended the stability analysis in this critically loaded regime to unbounded perturbations, establishing that stability can sometimes hold for the queue fed by scheduled arrivals, while the time-reversed scheduled arrival process can make the same queue unstable. In \citet{AramanGlynn12}, a functional central limit theorem (FCLT) for scheduled traffic is derived when the perturbations have infinite mean and Pareto-like tails. The limit involves a fractional Brownian motion (fBM) with Hurst parameter $H$ in $(0,1/2)$, from which a heavy traffic limit process for the workload can be obtained. \cite{AramanGlynn22} establish properties of scheduled traffic and show, for finite-mean Pareto-like perturbations, that an $S/D/1$ queue behaves very differently from both a $D/D/1$ and a $G/D/1$ queue.

All the above papers establish, in different ways, that scheduled traffic is more regular and less variable than that associated with conventional traffic models (e.g. Markov modulated Poisson processes, renewal arrivals, etc). For example, the fBMs that arise in connection with scheduled traffic are processes that exhibit long-range negative correlations, and are more regular than is Brownian motion itself. In this paper, we show that when the perturbations are super-heavy-tailed (i.e. the logarithm of the perturbations has a Pareto-like tail), the negative correlations are completely eliminated, in the sense that the associated FCLT has Brownian motion (H = 1/2) as its limit process; see Theorem~1. This paper also shows that when the perturbations have a Cauchy tail, the associated FCLT yields a Gaussian limit process that can be interpreted as an fBM with $H = 0$; see Theorem 2. These theorems round out our understanding of scheduled traffic, and establish that when the perturbations have infinite mean, the entire range of fBMs with $H$ in  $[0,1/2]$ are possible as limits, including the end points 0 and 1/2. 

\section{Description of Main Results}\label{Sec: Description}
Let $(\xi_j:j\in\mathbb{Z})$ be an i.i.d. sequence of perturbations. We view $\xi_j$ as the perturbation associated with the customer scheduled to arrive at time $j h$ (with $h>0$) so that its actual arrival time is $jh+\xi_j$. If $\tilde N_0$ is the random measure for which $$\tilde N_0(A)=\sum_jI(jh+\xi_j\in A),$$
for measurable $A\subseteq \mathbb R$, then $\tilde N_0(A)$ is the number of customers to arrive in the subset $A$. The random measure $\tilde N_0$ is $\mathbb{Z}$-stationary, in the sense that $\tilde N_0(\cdot+nh)\overset{\cal D}{=}\tilde N_0(\cdot)$ for $n\in\mathbb{Z}$ (where $A+t\overset{\Delta}{=}\{x+t:x\in A\}$, $\overset{\cal D}{=}$ denotes equality in distribution, and $\overset{\Delta}{=}$ denotes equality by definition). The scheduled traffic arrival counting process $N_0=(N_0(t):t\geq 0)$ is then defined via $N_0(t)\overset{\Delta}{=}\tilde N_0((0,t]).$

If we prefer a (fully) time-stationary version of the scheduled traffic process, we introduce a rv $U$ uniform on $[0,1]$ independent of $\xi_0,$ and define the random measure $\tilde N_1$ for which $$\tilde N_1(A)=\sum_jI(jh+Uh+\xi_j\in A).$$ Then, $\tilde N_1$ is time-stationary, in the sense that $\tilde N_1(\cdot+t)\overset{\cal D}{=}\tilde N_1(\cdot)$ for $t\in\mathbb{R}$. Its associated counting process is given by $N_1(t)\overset{\Delta}{=}\tilde N_1((0,t])$ for $t\geq 0$. In \citet{AramanGlynn22}, it is shown that \textit{regardless of the distribution of $\xi_0$}, the following properties hold
\begin{itemize}
\item [$i.)$] For each $t\geq 0$, $N_1(t)$ is light-tailed and, in fact $\e\exp(\theta N_1(t))<\infty$ for each $\theta\in\mathbb{R};$

\item [$ii.)$] $\e N_1(t)=t/h$ for $t\geq 0.$
\end{itemize}

Recall that the \textit{dispersion} of a non-negative rv Y is given by $\var Y/\e Y.$ Observe that 
\begin{align*}
    \var N_0(nh)&=\var (\sum_jI(jh+\xi_j\in (0,nh]))\\
    &= \sum_j\var I(jh+\xi_j\in (0,nh])\\
    &\leq  \sum_j\e I(jh+\xi_j\in (0,nh]))\\
    &=\e N_0(nh),
\end{align*}
so the dispersion of $N_0(nh)$ is always less than or equal to 1. This implies, in particular, that $\var N_0(nh)\leq \e N_0(nh)=n$, so the standard deviation of $N_0(nh)$ can never grow faster than $n^{1/2}$, regardless of the distribution of $\xi_0$. 

Let 
$$\cE(nh)=\sum_{j>n}I(jh+\xi_j\leq nh)$$

and

$$\cL(nh)=\sum_{j\leq n}I(jh+\xi_j> nh)$$ be the number of customers at time $nh$ that have arrived early (i.e. have been scheduled to arrive after $nh$ but have already arrived) or will arrive late (i.e. were scheduled to arrive before $nh$, but have not yet arrived). The sequence $((\cE(nh),\cL(nh)):n\in\mathbb{Z})$ is clearly $\mathbb{Z}$-stationary. But note that the Borel-Cantelli lemma and its converse imply that $\cE(nh)<\infty~a.s.$ if and only if $\e\xi_0^-<\infty,$ whereas $\cL(nh)<\infty~a.s.$ if and only if $\e \xi_0^+ <\infty$. (Of course, stationary versions of $\tilde N_0$ and $\tilde N_1$ exist regardless of whether $\e|\xi_0|$ is finite or not.)

When $\e|\xi_0|<\infty$, 
\begin{align}\label{Eq: N(nh)-nh}
  N_0(nh)-n&=\cE(nh)-\cL(nh)-(\cE(0)-\cL(0))\nonumber\\
  &\Rightarrow \cE'(0)-\cL'(0)-(\cE(0)-\cL(0))
    \end{align}
as $n\rightarrow\infty,$ where $((\cE'(0),\cL'(0)))\overset{\cal D}{=}((\cE(0),\cL(0)))$ and $((\cE'(0),\cL'(0)))$ is independent of $((\cE(0),\cL(0)))$ ; see Theorem~1 of \citet{AramanGlynn22} for the corresponding proof for $N_1(\cdot)$ (the proof for $N_0$ is easier).

Thus, when $\xi_0$ has finite mean, $N_0(nh)-n$ is stochastically bounded. We turn next to describing the behavior when $\xi_0$ has infinite mean. To simply describe the result, we henceforth assume that $h=1$ and that the $\xi_j$'s are positive rvs. Suppose that
\begin{equation}\label{Eq: Pareto_Assumption}
    \Pro(\xi_0>x)\sim \kappa\,x^{-r}
\end{equation}
as $x\rightarrow\infty$ for $0<\kappa<\infty$ and $0<r<1$. For $t\geq 0$, set 
$$\chi_n(t)=\frac{N_0(nt)-\lfloor n t\rfloor}{n^{(1-r)/2}}$$
where $\lfloor x \rfloor$ is the floor of $x$ (i.e. the greatest integer less than or equal to $x$). Suppose that $B_H=(B_H(t):t\geq 0)$ is a fBm with unit variance and Hurst parameter $H\in(0,1),$ so that it is a mean zero Gaussian process with covariance function given by 
\begin{equation}\label{Eq: Covariance}
\cov (B_H(s), B_H(t))=\frac{1}{2}(|s|^{2H}+|t|^{2H}-|t-s|^{2H})
\end{equation}
for $s,t\geq 0$. Put $\chi_n=(\chi_n(t):t\geq 0).$ In \citet{AramanGlynn12}, it is shown that 
$$\chi_n\Rightarrow \sqrt{2\kappa(1-r)^{-1}}\, B_H$$
as $n\rightarrow\infty$ in $D[0,\infty)$, where $H=(1-r)/2.$ Hence, as the perturbation tails get heavier, the fluctuations of $X_n$ become larger, approaching the $n^{1/2}$ theoretical limit imposed by the unit dispersion upper bound discussed earlier. Furthermore, the limit process $B_H$ becomes more ``disordered", losing the strong negative correlations that are present when $H$ is close to 0.

This raises the question of whether there exist perturbation distributions that achieve the disordered state that is associated with Brownian motion, namely $H=1/2$. We say that a positive rv $\xi$ has a \textit{super-heavy tail} if there exists a (deterministic) function $a(t)\leq t$ for which 
\begin{equation}\label{eq: super_heavy1}
    \frac{\Pro(\xi>a(t))}{\Pro(\xi>t)}\rightarrow 1
\end{equation}
as $t\rightarrow\infty,$ and 
\begin{equation}\label{eq: super_heavy2}
    \frac{a(t)}{t\,\Pro(\xi>a(t))}\rightarrow 0
\end{equation}
as $t\rightarrow\infty$. Put $\bar F(t)\overset{\Delta}{=}\Pro(\xi_0>t)$. We note that if $\bar F(t)\sim c\,(\log t)^{-\alpha}$ as $t\rightarrow\infty$ for $c,\alpha>0,$ then $a(t)=t^b$ for $b\in(0,1)$ satisfies (\ref{eq: super_heavy1}) and (\ref{eq: super_heavy2})  and hence $\xi_0$ is super-heavy tailed. It is also evident that  (\ref{eq: super_heavy1}) requires that $a(t)\rightarrow\infty$ as $t\rightarrow\infty$, while  (\ref{eq: super_heavy2}) implies that $a(t)/t\rightarrow 0$ as $t\rightarrow\infty$.

Let $B=(B(t):t\geq 0)$ be standard Brownian motion. For $t\geq 0$, put $$X_n(t)=\frac{1}{\sqrt{n\,\bar F (n)}}(N_0(nt)-\lfloor n t\rfloor).$$
and set $X_n=(X_n(t):t\geq 0).$
\begin{theorem}\label{th:Conv_to_Bm}
    Suppose that $\xi_0$ is a positive rv that is super-heavy tailed. Then  $$X_n\Rightarrow \sqrt{2}B$$ as $n\rightarrow\infty$ in $D[0,\infty).$
\end{theorem}

The proof is deferred to Section~\ref{sec:Section_3}. Theorem~\ref{th:Conv_to_Bm} shows that when $\xi_0$ is super-heavy tailed, then the limiting behavior of $X_n$ is Brownian, and the stochastic fluctuations of $N_0(n)$ can come arbitrarily close to the $n^{1/2}$ theoretical limit. (We can take, for example, $\bar F(t)=(\log_k t)^\alpha,$ where $\log_1 t=\log t$, and $\log_k t=\log(\log_{k-1} t)$ for $k\geq 1$.) So $H=1/2$ is achievable as a Hurst parameter arising from scheduled traffic.  Note that $N_1(\cdot)\overset{\cal D}{=}N_0(1-U+\cdot)-N_0(1-U),$ so that $N_1(\cdot)$ satisfies the same FCLT as does $N_0(\cdot)$ in Theorem~\ref{th:Conv_to_Bm}.

We also take this opportunity to explore further the behavior of scheduled traffic as one transitions from the finite mean setting to the context of (\ref{Eq: Pareto_Assumption}). In particularly, suppose that $\xi_0$ has Cauchy-like tails, so that 
\begin{equation}\label{eq: Cauchy_dist}
    \Pro(\xi_0>x)\sim d/x
\end{equation}
as $x\rightarrow\infty$, when $d>0$. 
Put $$X'_n(t)=\frac{N_0(nt)-\lfloor n t\rfloor}{\sqrt{\log n}}$$ for $t\geq 0$. Let $Z=(Z(t):t\geq 0)$ be a zero mean Gaussian process with $\var Z(t)=2$ for $t>0$ with $\cov(Z(s),Z(t))=1 $ for $s\neq t.$ We write $X'_n\overset{fdd}{\rightarrow}Z$ as $n\rightarrow\infty$ if the finite-dimensional distributions of $X'_n$ converge weakly to those of $Z$ as $n\rightarrow\infty.$ 
\begin{theorem}\label{th:conv_to_std_normal}
    If (\ref{eq: Cauchy_dist}) holds, then 
    \begin{equation}\label{eq: Cauchy_limit}
    X'_n(\cdot)\overset{fdd}{\rightarrow}\sqrt{d}\,Z(\cdot)
    \end{equation}
as $n\rightarrow\infty$.
\end{theorem}
In this Cauchy setting, the stochastic fluctuations of $N_0(n)$ are small (of order $(\log n)^{1/2}$) and the limit process $Z$ exhibits such strong negative autocorrelations that $\var Z(t)$ does not grow as $t\rightarrow\infty.$ Note that $Z(t)-Z(s)\overset{\cal D}{=}Z(t-s)\overset{\cal D}{=}{\cal N}(0,2),$ regardless of how close $s$ is to $t$. Hence, Lemma~7.7, p. 131, of \citet{Ethier86} implies that $Z\notin D[0,\infty)$ so $(X'_n:n\geq 1)$ is not tight in $D[0,\infty).$

Although there does not appear to be a definition of fractional Brownian motion at $H=0$ given by the literature, the limit appearing in (\ref{eq: Cauchy_limit}) is consistent with the covariance function (\ref{Eq: Covariance}) when $H=0$. We may therefore (reasonably) take the view that Theorem~\ref{th:conv_to_std_normal}'s limit process corresponds to a fractional Brownian motion with $H=0$.  

\section{Proof of Main Results}\label{sec:Section_3}
\subsection{Proof of Theorem~\ref{th:Conv_to_Bm}}
Put $a_n=(\bar F(n)\,n)^{-1/2}$. For $\theta_1,\theta_2\in \mathbb{R}$, put $\tilde \theta_1=a_n\theta_1$ and $\tilde \theta_2=a_n\theta_2$. Let 
$n_1=\lfloor n s \rfloor$, $n_2=\lfloor n t\rfloor,$ for $0\leq s\leq t$. Note that 
\begin{align*}
X_n(t) &=  a_n\big(\sum_{j=1}^{\lfloor nt\rfloor}I(j + \xi_j \in
(0,nt])-\lfloor nt\rfloor\nonumber\\
& ~~~~~~~~~~~~~~~~~~~~ + \sum_{j\leq 0} I(j+\xi_j \in (0,nt])\big) \nonumber\\
&=  a_n\big(-\sum_{j=1}^{\lfloor nt\rfloor}I(j + \xi_j >nt) \nonumber\\
& ~~~~~~~~~~~~~~~~~~~~ + \sum_{j\leq 0} I(j+\xi_j \in (0,nt])\big).\nonumber
\end{align*}
Hence,
\begin{align*}
&~\theta_1 X_n(s)+\theta_2 X_n(t)\nonumber\\
=&~-(\tilde{\theta}_1+\tilde{\theta}_2)\sum_{j=1}^{n_1}I(\xi_j +j>n_2)-\tilde{\theta}_1\sum_{j=1}^{n_1}I(\xi_j+j\in (n_1,n_2])-\tilde{\theta}_2\sum_{j= n_1 +1}^{n_2}I(\xi_j +j>n_2) \nonumber\\
&~+(\tilde{\theta}_1+\tilde{\theta}_2)\sum_{j\leq 0}I(\xi_j +j\in
(0,n_1])+\tilde{\theta}_2\sum_{j\leq 0}I(\xi_j +j\in (n_1,n_2]).\nonumber
\end{align*}

Then, the joint log moment generating of $(X_n(s),X_n(t))$ evaluated at $(\theta_1,\theta_2)$ is given by 
\begin{align*}
&\sum_{j=1}^{ n_1
}\log\left(1+(e^{-\tilde{\theta}_1-\tilde{\theta}_2}-1)\bar{F}(n_2-j)
+(e^{-\tilde{\theta}_1}-1)(\bar{F}(n_1-j)-\bar{F}(n_2-j))\right) \nonumber\\
+&\sum_{j= n_1  +1}^{
n_2  }\log\left(1+(e^{-\tilde{\theta}_2}-1)\bar{F}(n_2-j)\right) \nonumber\\
+&\sum_{j\leq0}\log\left(1+(e^{\tilde{\theta}_1+\tilde{\theta}_2}-1)(\bar{F}(-j)-\bar{F}(n_1-j))+(e^{\tilde{\theta}_2}-1)(\bar{F}(n_1-j)-\bar{F}(n_2-j))\right).
\end{align*}\label{Equ:LogMmtGenfct}
Note that $\tilde{\theta_i}\rightarrow0$ as $n\rightarrow\infty$ and $\log(1+x)=x(1+o(1))$ as $x\rightarrow0.$ Hence
\begin{align*}
&~\log \mathbb{E}\exp (\theta_1 X_n(t_1)+\theta_2 X_n(t_2))\nonumber\\
=&~(e^{-\tilde{\theta}_1-\tilde{\theta}_2}-e^{-\tilde{\theta}_1})\sum_{j=1}^{ n_1}\bar{F}(n_2-j)
+(e^{-\tilde{\theta}_1}-1)\sum_{j=1}^{ n_1}\bar{F}(n_1-j) +(e^{-\tilde{\theta}_2}-1)\sum_{j= n_1+1}^{n_2}\bar{F}(n_2-j) \nonumber\\
&~~~~~+(e^{\tilde{\theta}_1+\tilde{\theta}_2}-1)\sum_{j=0}^{\infty}(\bar{F}(j)-\bar{F}(n_1+j))+ (e^{\tilde{\theta}_2}-1)\sum_{j=0}^{\infty}(\bar{F}(n_1+j)-\bar{F}(n_2+j))+o(1)
\end{align*}
as $n\rightarrow\infty.$
For $0\leq k_1\leq k_2$ and $r\geq k_2-k_1$,
\begin{align}
\sum_{j=0}^r(\bar{F}(k_1+j)-\bar{F}(k_2+j))
=\sum_{j=k_1}^{k_2-1}\bar{F}(j)-\sum_{j=k_1+r+1}^{k_2+r}\bar{F}(j).\nonumber
\end{align}
Because, $\sum_{j=k_1+r+1}^{k_2+r}\bar{F}(j)\leq (k_2-k_1)\bar{F}(k_1+r+1)\rightarrow0$ as $r\rightarrow\infty$, it follows that
\begin{equation}
\sum_{j=0}^\infty(\bar{F}(k_1+j)-\bar{F}(k_2+j))=\sum_{j=k_1}^{k_2-1}\bar{F}(j).\nonumber
\end{equation}
even when $\sum_{j=k_1}^\infty\bar F(j)=\infty.$ Also, $a_n\,(\e N(nt)-\lfloor nt\rfloor)=o(1)$,  when $n\rightarrow\infty$ uniformly in $t\geq 0$, so 

\begin{align}
&~\log \mathbb{E}\exp (\theta_1 X_n(t_1)+\theta_2 X_n(t_2))\nonumber\\
=&~(e^{-\tilde{\theta}_1-\tilde{\theta}_2}-e^{-\tilde{\theta}_1})\sum_{j=n_2-n_1}^{ n_2-1}\bar{F}(j)
+(e^{-\tilde{\theta}_1}-1)\sum_{j=0}^{ n_1-1}\bar{F}(j) \nonumber\\
&~~~~~+(e^{-\tilde{\theta}_2}-1)\sum_{j=0}^{n_2-n_1-1}\bar{F}(j)+(e^{\tilde{\theta}_1+\tilde{\theta}_2}-1)\sum_{j=0}^{n_1-1}\bar{F}(j)\nonumber\\
&~~~~~+ (e^{\tilde{\theta}_2}-1)\sum_{j=n_1}^{n_2-1}\bar{F}(j)+o(1)\nonumber\\
=&~\frac{1}{2}((\theta_1+\theta_2)^2-\theta_1^2))\,a^2_n\sum_{j=n_2-n_1}^{ n_2-1}\bar{F}(j)\nonumber\\
&~~~~~+\frac{1}{2}\theta_1^2\,a^2_n\sum_{j=0}^{ n_1-1}\bar{F}(j) +\frac{1}{2}\theta_2^2a^2_n\sum_{j=0}^{n_2-n_1-1}\bar{F}(j)\nonumber\\
&~~~~~+\frac{1}{2}(\theta_1+\theta_2)^2)a^2_n\sum_{j=0}^{n_1-1}\bar{F}(j)\nonumber\\
&~~~~~+ \frac{1}{2}\,\theta_2^2\,a^2_n\sum_{j=n_1}^{n_2-1}\bar{F}(j)+o(1)\label{eq:logmomtgen_a_n}
\end{align}
as $n\rightarrow\infty.$ For $r>0$, note that 
$$a_n^2\sum_{j=0}^{\lfloor rn\rfloor}\bar F(j)=\frac{1}{n\bar F(n)} \sum_{j=0}^{\lfloor a(rn)\rfloor}\bar F(j)+\frac{1}{n}\sum_{\lceil a(rn)\rceil}^{\lfloor rn\rfloor}\frac{\bar F(j)}{\bar F(n)}.$$
Because $a(n)/(r n)\rightarrow 0$ as $n\rightarrow\infty,$ 
\begin{align*}
    \frac{1}{n\bar F(n)} \sum_{j=0}^{\lfloor a(rn)\rfloor}\bar F(j)&\leq \frac{a(rn)}{n\,\bar F(n)}=\frac{a(rn)}{rn\bar F(rn)}\cdot\frac{r\,\bar F(rn)}{\bar F(n)}\\
    &\leq \frac{a(rn)}{rn\bar F(rn)}\cdot\frac{r\,\bar F(a(n))}{\bar F(n)}\rightarrow 0\cdot r=0
\end{align*}
as $n\rightarrow\infty$. Again, because $a(rn)/n\rightarrow 0,$ it follows that if $r>1$, 
$$\bar F(a(rn))/\bar F(rn)\geq \bar F(n)/\bar F(rn)\geq 1,$$ so $\bar F(n)/\bar F(rn)\rightarrow 1$ as $n\rightarrow\infty.$ Similarly, $\bar F(n)/\bar F(rn)\rightarrow 1$ if $0<r\leq 1.$ Consequently, 
$$1\leq \max_{a(rn)\leq t\leq n}\frac{\bar F(t)}{\bar F(n)}\leq \frac{\bar F(a(rn))}{\bar F(rn)}\cdot\frac{\bar F(rn)}{\bar F(n)}\rightarrow 1\cdot 1=1, $$
and 
$$\left(\frac{\lfloor rn\rfloor-\lceil a(rn)\rceil }{n}\right)\leq \frac{1}{n}\sum_{j=\lceil a(rn)\rceil }^{\lfloor rn\rfloor}\frac{\bar F(j)}{\bar F(n)}\leq r\max_{a(rn)\leq t\leq n}\frac{\bar F(t)}{\bar F(n)}.$$
It follows that $$\frac{1}{n}\sum_{j=\lceil a(rn)\rceil }^{\lfloor rn\rfloor}\frac{\bar F(j)}{\bar F(n)}\rightarrow r$$
as $n\rightarrow\infty$, so that $$a_n^2\sum_{j=0}^{\lceil a(rn)\rceil}\bar F(j)\rightarrow r$$ 
as $n\rightarrow\infty$. 

Hence, (\ref{eq:logmomtgen_a_n}) implies that the log moment generating function of $(X_n(s),X_n(t))$  converges to
\begin{align*}
    &~~~~\frac{1}{2}((\theta_1+\theta_2)^2-\theta_1^2))\,s+\theta_1^2/2\,s +\theta_2^2/2\,(t-s)\\
    &~~~+\frac{1}{2}(\theta_1+\theta_2)^2 \,s+\theta_2^2/2\,(t-s)\\
    &=(\theta_1+\theta_2)^2\,s+\theta_2^2\,(t-s),
\end{align*}
which is the log moment generating function of $\sqrt{2}\,(B(s),B(t)).$ A similar  calculation establishes that the finite-dimensional distribution of $X_n$ converge to those of $\sqrt{2} \,B.$

We note that $N_0(n)$ is the sum of the $n$ stationary increments $\tilde N_0((i-1,i])$ ($1\leq i\leq n$). Furthermore, because we have established convergence of the moment generating function for the finite-dimensional distributions, it follows that $$\e \left(\frac{N_0(n)^2}{n\bar F(n)}\right)^k\rightarrow 2^k$$
as $n\rightarrow\infty$ for each $k\geq 1$. As argued above, $\bar F(r t)/\bar F(t)\rightarrow 1$ as $t\rightarrow \infty$ for $r>0,$ so $\bar F$ is slowly varying. 
It follows from Lemma~2.1 of \citet{Taqqu75} that 
$$\left(\frac{N_0(\lfloor n s\rfloor)-\lfloor n s\rfloor)}{\sqrt{n\,\bar F (n)}}:0\leq s\leq t\right)$$ is tight in $D[0,t]$ for each $t\geq 0$. Since 
\begin{equation}\label{eq: tight}
    \max_{0\leq s\leq t}\frac{N_0(n s)-N_0(\lfloor n s\rfloor)}{\sqrt{n\,\bar F (n)}}\leq \max_{1\leq i\leq \lceil nt\rceil}\frac{N_0(i)-N_0(i-1)}{\sqrt{n\,\bar F (n)}},
\end{equation}
the stationarity of $(N_0(i)-N_0(i-1):i\geq 1)$ and the finiteness  of their common moment generating function imply that the right-hand side of (\ref{eq: tight}) converges to 0 almost surely. Hence, $(X_n:n\geq 1)$ is tight in $D[0,t]$ for each $t\geq 0$, proving the theorem. \qedsymbol{}

\subsection{Proof of Theorem~\ref{th:conv_to_std_normal}}
Put $a_n=(\log n)^{-1/2}$. We follow the identical steps of the proof of Theorem~\ref{th:Conv_to_Bm} up to (\ref{eq:logmomtgen_a_n}). Now recall that $\sum_{j=1}^n 1/j=\log n+\gamma+o(1)$ as $n\rightarrow\infty$ (see, p. 388, of \citet{widder89}).

Then for $r>0$,
\begin{align*}
a_n^2\sum_{j=0}^{\lfloor rn\rfloor}\bar F(j)&=\frac{1}{\log n}\sum_{j=0}^{\lfloor a_n^{-1}\rfloor}\bar F(j)+\frac{1}{\log n}\sum_{j=\lfloor a_n^{-1}\rfloor+1}^{\lfloor rn\rfloor}\bar F(j)\\
&= O(\frac{1}{\sqrt{\log n}})+\frac{1}{\log n}\sum_{j=\lfloor a_n^{-1}\rfloor+1}^{\lfloor rn\rfloor}\frac{d}{j} (1+o(1))\\
&=O(\frac{1}{\sqrt{\log n}})+\frac{1}{\log n}\left(\sum_{j=1}^{\lfloor rn\rfloor}\frac{d}{j} -\sum_{j=1}^{\lfloor a_n^{-1}\rfloor}\frac{d}{j}\right)(1+o(1))\\
&=O(\frac{1}{\sqrt{\log n}})+\frac{d}{\log n}\Big(\log (rn)-1/2 \log\log n\Big)(1+o(1))\\
&=d+o(1)
\end{align*}
as $n\rightarrow\infty$. 
It follows that the log moment generating function of $(X'_n(s),X'_n(t))$ evaluated at $(\theta_1,\theta_2)$ converges to 
\begin{align*}
    &[\theta_1^2/2+(\theta_1-\theta_2)^2/2+\theta_2^2/2]\\
    &=\frac{d}{2}\,[2\theta_1^2-2 \theta_1\theta_2+2\theta_2^2],
\end{align*}
which is the log moment generating function of $(Z(s),Z(t))$. A similar argument establishes convergence of the finite-dimensional distributions in the general case.
\qedsymbol{}

\bibliographystyle{ormsv080}
\bibliography{BiblioApP}
\end{document}